\documentclass[preprint,3p,sort&compress,final,times]{elsarticle}
\usepackage{amsmath}
\usepackage{amssymb}
\usepackage{amsthm}
\usepackage{mathtools}
\usepackage{placeins}
\usepackage{bm}
\usepackage{caption}
\usepackage{subcaption}
\usepackage[usenames,dvipsnames]{color}

\usepackage{cases}

\usepackage{lineno}

\usepackage{tikz-cd}

\newcommand{\tabref}[1]{Table~\ref{#1}} 
\newcommand{\figref}[1]{Figure~\ref{#1}} 

\newcommand{\transpose}{\intercal}                      

\newcommand{\bb}[1]{\left( #1 \right)}

\newcommand{\nl}{\\[1.1ex]}

\newcommand{\VJ}[1]{\textcolor{black}{#1}}
\theoremstyle{plain}

\journal{}
\begin{document}
\begin{frontmatter}

\title{Estimation and numerical validation of inf-sup constant for \VJ{bilinear form $(p, \mathrm{div}\ \bm{u})$}}
\author[TUD]{V.~Jain\corref{cor}}
\ead{V.Jain@tudelft.nl}
\author[TUD]{M.~Gerritsma}

\address[TUD]{Delft University of Technology, Faculty of Aerospace Engineering, P.O. Box 5058, 2600 GB Delft, The Netherlands}

\cortext[cor]{Corresponding author. Tel. +31 15 2789670.}

\begin{abstract}
We give a derivation for the value of inf-sup constant for the bilinear form $\bb{p \;, \mathrm{div}\ \bm{u}}$.
We prove that the value of inf-sup constant is equal to 1.0 in all cases and is independent of the size and shape of the domain.
Numerical tests for validation of inf-sup constants is performed using finite dimensional spaces defined in \cite{2020jain} on two test domains i) a square of size $\Omega = [0,1]^2$, ii) a square of size $\Omega = [0,2]^2$, for varying mesh sizes and polynomial degrees.
The numeric values are in agreement with the theoretical value of inf-sup term.
\end{abstract}


\end{frontmatter}
\section{Introduction}
\VJ{Numerical schemes for mixed finite element methods} should result in stable bounded solutions.
The stability of finite element discretization is governed by the inf-sup criterion for bilinear forms.
The objective of this paper is to derive and validate the inf-sup stability constant for the bilinear form of the divergence term.
This is of importance because it often appears as \VJ{pressure constraint} term in discretization of fluid mechanics equations.
For given $\bm{u} \in H\bb{\mathrm{div}; \Omega}$, $p \in L^2\bb{\Omega}$, the bilinear form is give by
\begin{equation} \label{eq:bilinear_form_div}
b\bb{p,\bm{u}} = \bb{p \;, \mathrm{div}\ \bm{u}} \;.
\end{equation}
The inf-sup constant $\beta$ for \eqref{eq:bilinear_form_div} is given by
\begin{equation}
\beta = \inf_{p \in L^2\bb{\Omega} } \sup_{\bm{u} \in H\bb{\mathrm{div};\Omega}} \frac{\bb{p \;, \mathrm{div}\ \bm{u}}}{\| p \|_{L^2\bb{\Omega}} \|\bm{u} \|_{H\bb{\mathrm{div};\Omega}} } \;.
\end{equation}
\VJ{Quite some work has been done to approximate the value of inf-sup constant see for eg. \cite{chapelle1993inf, bathe2001inf, bernardi2016continuity}.}
\VJ{In this work we measure the norm of the velocity field in the space orthogonal to kernel of the divergence operator} and prove that the value of the inf-sup constant $\beta = 1.0$.
For numerical tests we use finite dimensional spaces defined in \cite{2020jain}.
The validation study is performed on three different domains: i) unit square $\Omega = [0,1]^2$; ii) square $\Omega = [0,2]^2$, for varying mesh sizes, $h$, and polynomial degrees $N=1,2,3$.
\VJ{The numerical values are in agreement with the derived values of $\beta$, with maximum errors shown of the order of $10^{-6}$. }

\section{Derivation of inf-sup constant $\beta$}
\VJ{Let $\Omega \subset \mathbb{R}^2$ be an open, bounded domain.}
We will use the finite dimensional spaces \VJ{and the divergence operator} defined in \cite[\S3]{2020jain}: $D\bb{\Omega} \subset H\bb{\mathrm{div};\Omega}$, $\widetilde{S}\bb{\Omega} \subset L^2\bb{\Omega}$, $S\bb{\Omega} \subset L^2\bb{\Omega}$, \VJ{and $\mathbb{E}^{2,1}$ the discrete representation of the divergence operator}.
Let $K = \mathrm{Ker}\ \mathbb{E}^{2,1}$, $H = \mathrm{Ker}\bb{{\mathbb{E}^{2,1}}^\transpose} $.
The discrete inf-sup condition is then given by
	\begin{equation} \label{eq:discrete_inf_sup}
	\beta = \inf_{p \in H^\perp} \sup_{\bm{u} \in K^\perp} \frac{\mathcal{N}^1\bb{\bm{u}} ^\transpose {\mathbb{E}^{2,1}}^\top \widetilde{\mathcal{N}}^0\bb{p}}{\| \bm{u} \|_{H(\mathrm{div};\Omega)} \|p\|_{L^2\bb{\Omega}}} \;,
	\end{equation}
	\VJ{where $\mathcal{N}^x$ denotes the degrees of freedom in our finite element, i.e. the vector of expansion coefficients}.
	Here $p \in \widetilde{S}(\Omega) \subset L^2(\Omega)$ and $\bm{u} \in D(\Omega) \subset  H(\mathrm{div};\Omega)$ and the vectors $\widetilde{\mathcal{N}}^0 \bb{p}$ and $\mathcal{N}^1\bb{\bm{u}}$ are the expansion coefficients.
The norm of $p \in \widetilde{S}(\Omega)$ is
\begin{equation} \label{eq:norm_p1}
 \| p \|_{L^2\bb{\Omega}}^2 = {\widetilde{\mathcal{N}}^0\bb{p}}^\transpose { \mathbb{M}^{(2)} }^{-1} \widetilde{\mathcal{N}}^0\bb{p} \;,
 \end{equation}
\VJ{ the norm of $\bm{u} \in H\bb{\mathrm{div};\Omega}$ is
 \[ \| \bm{u} \|^2_{H\bb{\mathrm{div};\Omega}} =  {\mathcal{N}^1\bb{\bm{u}}} ^\transpose \bb{\mathbb{M}^{(1)} + {\mathbb{E}^{2,1}}^\transpose \mathbb{M}^{(2)} \mathbb{E}^{2,1} } \mathcal{N}^1\bb{\bm{u}} \;, \]}
	and the norm of $\bm{u} \in K^\perp \subset D(\Omega)$ is
	\begin{equation} \label{eq:norm_u1}
\| \bm{u} \|_{K^\perp}^2 = \mathcal{N}^1\bb{\bm{u}}^\transpose {\mathbb{E}^{2,1}}^\transpose \mathbb{M}^{(2)} {\mathbb{E}^{2,1}} \mathcal{N}^1\bb{\bm{u}} = \| \mathrm{div}\ \bm{u} \|_{L^2\bb{\Omega}} \;.
	\end{equation}
	In the continuous case, using Cauchy Schwartz inequality, we have, for all $p \in L^2\bb{\Omega}$ and $\bm{u} \in K^\perp$
	\begin{equation} \label{eq:ineq1}
\VJ{	\frac{\left (p, \mathrm{div}\ \bm{u} \right )}{\| \bm{u} \|_{K^\perp} \|p\|_{L^2\bb{\Omega}}} \leq \frac{\| \mathrm{div}\ \bm{u})\|_{L^2\bb{\Omega}} \|p\|_{L^2\bb{\Omega}}}{\| \bm{u} \|_{K^\perp} \|p\|_{L^2\bb{\Omega}}} = 1 \;. }
	\end{equation}
If this inequality holds for all $p$ and $\bm{u}$ it should also hold when we take the infimum over $\widetilde{S}(\Omega)$ and the supremum over $D(\Omega)$, from which we conclude that $\beta \leq 1$.
	
Now, for an arbitrary vector field $\bm{u}^*$, we have
	\begin{equation}
	\sup_{u \in K^\perp} \frac{\mathcal{N}^1\bb{\bm{u}}^\transpose {\mathbb{E}^{2,1}}^\transpose \widetilde{\mathcal{N}}^0\bb{p}}{\| \bm{u} \|_{H(\mathrm{div};\Omega)} \|p\|_{L^2\bb{\Omega}}} \geq \frac{\mathcal{N}^1\bb{\bm{u}^*}^{\transpose} {\mathbb{E}^{2,1}}^\transpose \widetilde{\mathcal{N}}^0\bb{p} }{\| \bm{u}^* \|_{H(\mathrm{div};\Omega)} \|p\|_{L^2\bb{\Omega}}} \;,
	\label{eq:lower-bound}
\end{equation}
If we now take for $\bm{u}^*$, the vector field with expansion coefficients which satisfy $\mathbb{M}^{(2)} \mathbb{E}^{2,1} \mathcal{N}^1\bb{\bm{u}^*} = \widetilde{\mathcal{N}}^0\bb{p}$, then $\| p \|_{L^2\bb{\Omega}} = \| \bm{u}^* \|_{H(\mathrm{div};\Omega)}$ and the numerator $\mathcal{N}^1\bb{\bm{u}}^\transpose \mathbb{E}^\transpose \widetilde{\mathcal{N}}^0\bb{p} = \|p\|_{L^2\bb{\Omega}}^2 = \|\bm{u}^*\|_{H(\mathrm{div};\Omega)} \|p\|_{L^2\bb{\Omega}}$.
If we insert these estimates in \eqref{eq:lower-bound} we have
		\begin{equation} \label{eq:ineq2}
	\sup_{\bm{u} \in K^\perp} \frac{\mathcal{N}^1\bb{\bm{u}}^\transpose {\mathbb{E}^{2,1}}^\transpose \widetilde{\mathcal{N}}^0\bb{p} }{\| \bm{u} \|_{H(\mathrm{div};\Omega)} \|p\|_{L^2\bb{\Omega}}} \geq \frac{\mathcal{N}^1\bb{\bm{u}}^{*\transpose} {\mathbb{E}^{2,1}}^\transpose \widetilde{\mathcal{N}}^0\bb{p}}{\| \bm{u}^* \|_{H(\mathrm{div};\Omega)} \|p\|_{L^2\bb{\Omega}}} = \frac{\|\bm{u}^*\|_{H(\mathrm{div};\Omega)} \|p\|_{L^2\bb{\Omega}}}{\| \bm{u}^* \|_{H(\mathrm{div};\Omega)} \|p\|_{L^2\bb{\Omega}}} = 1\;,
	\end{equation}
which shows that $\beta \geq 1$.
From \eqref{eq:ineq1} and \eqref{eq:ineq2} we conclude that $\beta =1$ in the discrete setting. 
This value is independent of the mesh size or polynomial degree, so this value is also the inf-sup constant for $h \rightarrow 0$. This value is also independent of the domain $\Omega$.
\section{Evaluation of numeric inf-sup constant ${\beta}_h$}
In this section we will follow \cite[\S3.4.3]{2010Boffi} to evaluate the inf-sup constant.
Let $S_x$ and $S_y$ be the two symmetric and (semi-)positive definite matrices, such that
\begin{eqnarray}
S_x^\transpose S_x & = & {\mathbb{E}^{2,1}}^\transpose \mathbb{M}^{(2)} \mathbb{E}^{2,1} \;, \nl
S_y^\transpose S_y & = & {\mathbb{M}^{(2)}}^{-1} \;.
\end{eqnarray}
We can write the norms in \eqref{eq:norm_p1} and \eqref{eq:norm_u1} as
\begin{eqnarray}
|| \bm{u} || _{K^\perp} & = & {\mathcal{N}^1\bb{\bm{u}}}^\transpose {\mathbb{E}^{2,1}}^\transpose \mathbb{M}^{(2)} \mathbb{E}^{2,1} \mathcal{N}^1\bb{\bm{u}} = {\mathcal{N}^1\bb{\bm{u}}}^\transpose S_x^\transpose S_x \mathcal{N}^1\bb{\bm{u}} = || S_x \mathcal{N}^1\bb{\bm{u}} || _E \;, \nl
|| p ||_{L^2\bb{\Omega}} & = & {\widetilde{\mathcal{N}}^0\bb{p}}^\transpose {\mathbb{M}^{(2)}}^{-1} \widetilde{\mathcal{N}}^0\bb{p} = {\widetilde{\mathcal{N}}^0\bb{p}}^\transpose S_y^\transpose S_y \widetilde{\mathcal{N}}^0\bb{p} = ||S_y \widetilde{\mathcal{N}}^0\bb{p} ||_E \;,
\end{eqnarray}
where $|| \cdot ||_E$ is the Eucledian vector norm.

Let $M = S_y^{-1}\ \mathbb{E}^{2,1}\ S_x^{-1}$ and its singular value decomposition be given by
\begin{equation}  \label{eq:sigma}
M = S_y^{-1}\ \mathbb{E}^{2,1}\ S_x^{-1} = V \Sigma U \;.
\end{equation}
Now we can write $\mathbb{E}^{2,1}$ as
\begin{equation}\label{eq:transformedE}
\mathbb{E}^{2,1} = S_y S_y^{-1}\ \mathbb{E}^{2,1}\ S_x^{-1} S_x = S_y M S_x = S_y V \Sigma U S_x \;.
\end{equation}
Also, let
\begin{equation} \label{eq:u_and_p}
\mathcal{N}^1\bb{\bm{u}} = S_x^{-1} U^\transpose \bm{x} \;, \qquad \widetilde{\mathcal{N}}^0\bb{p} = S_y^{-1} V y \;.
\end{equation}
Now, if we substitute $\bm{u}$, $p$ from \eqref{eq:u_and_p} and $\mathbb{E}^{2,1}$ from \eqref{eq:transformedE} in the RHS term of \eqref{eq:discrete_inf_sup}, we get
\begin{eqnarray}
\inf_{p \in H^\perp} \sup_{\bm{u} \in K^\perp} \frac{{\widetilde{\mathcal{N}}^0\bb{p}} ^\transpose \mathbb{E}^{2,1} \mathcal{N}^1\bb{\bm{u}}}{||\bm{u}||_{K^\perp} ||p||_{L^2\bb{\Omega}} } & = & \inf_{y \in \bb{\mathrm{Ker}\ \Sigma^\transpose }^\perp} \sup_{\bm{x} \in \bb{\mathrm{Ker}\ \Sigma }^\perp} \frac{y^\transpose V^\transpose S_y^{-1} S_y V \Sigma U S_x S_x^{-1} U^\transpose \bm{x}}{|| S_x {S_x}^{-1} U^\transpose \bm{x} || _E || S_y {S_y}^{-1} V y ||_E} \nonumber \nl
& = & \inf_{y \in \bb{\mathrm{Ker}\ \Sigma^\transpose }^\perp} \sup_{\bm{x} \in \bb{\mathrm{Ker}\ \Sigma }^\perp} \frac{y^\transpose \Sigma \bm{x}}{|| \bm{x} || _E || y ||_E} := \beta _h \;. \label{eq:smallest_singular_value}
\end{eqnarray}
where in the last step we used \cite[Prop 3.4.3]{2010Boffi} which states that there exists a positive constant ${\beta}_h$ that is equivalent to the smallest positive singular value of the matrix, $M = S_y^{-1} E^{2,1} S_x ^{-1}$.
\begin{figure} [hbt]
\centering
\includegraphics[scale=0.4]{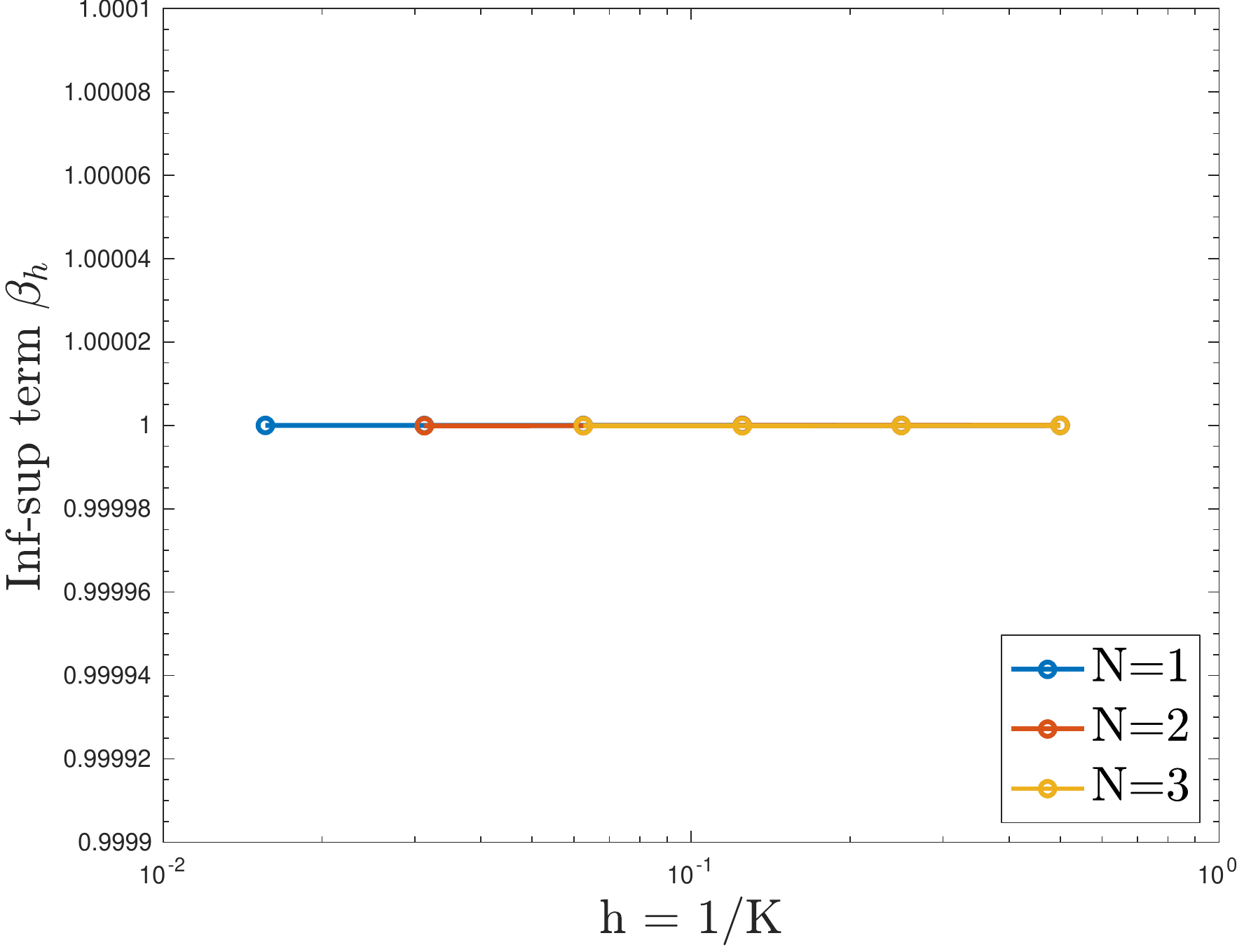}
\includegraphics[scale=0.4]{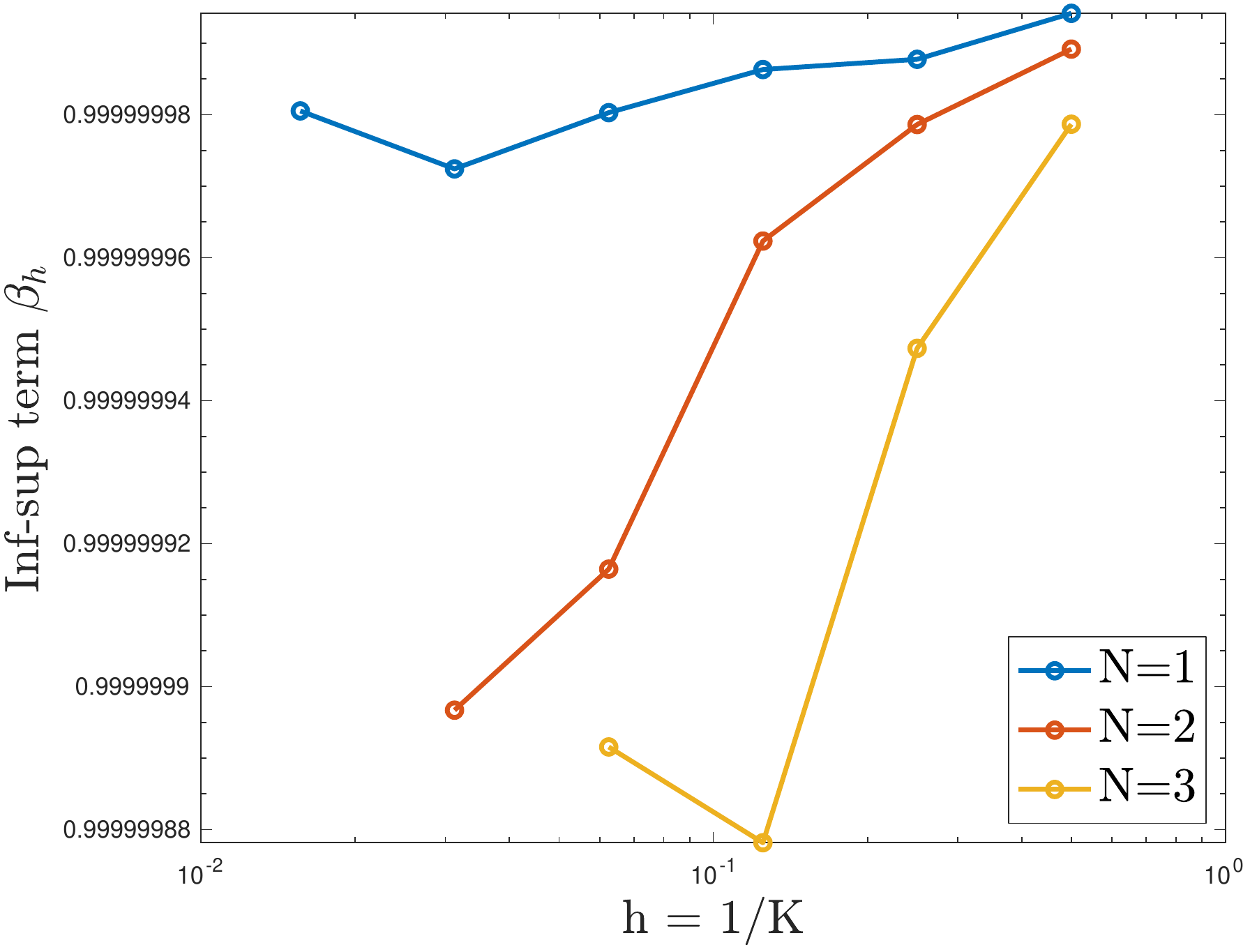}
\caption{Inf-sup term ${\beta}_h$ for $\Omega = [0,1]^2$ and polynomial degree $N=1,2,3$, with varying mesh refinement. Both the plots have same values. The plot on the right side is a zoom in of data.}
\label{fig:inf_sup_beta1}
\end{figure}
\begin{figure} [hbt]
\centering
\includegraphics[scale=0.4]{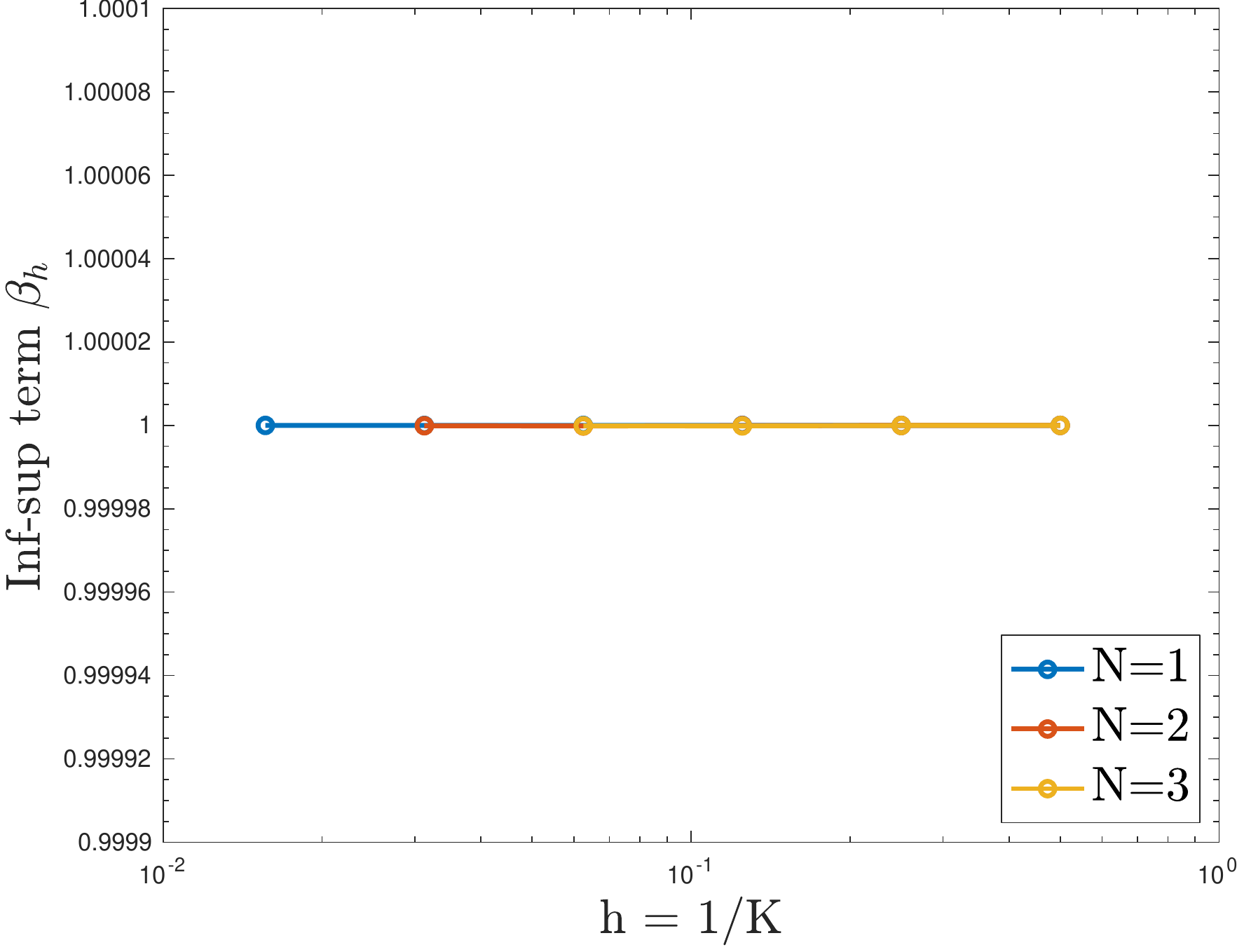}
\includegraphics[scale=0.4]{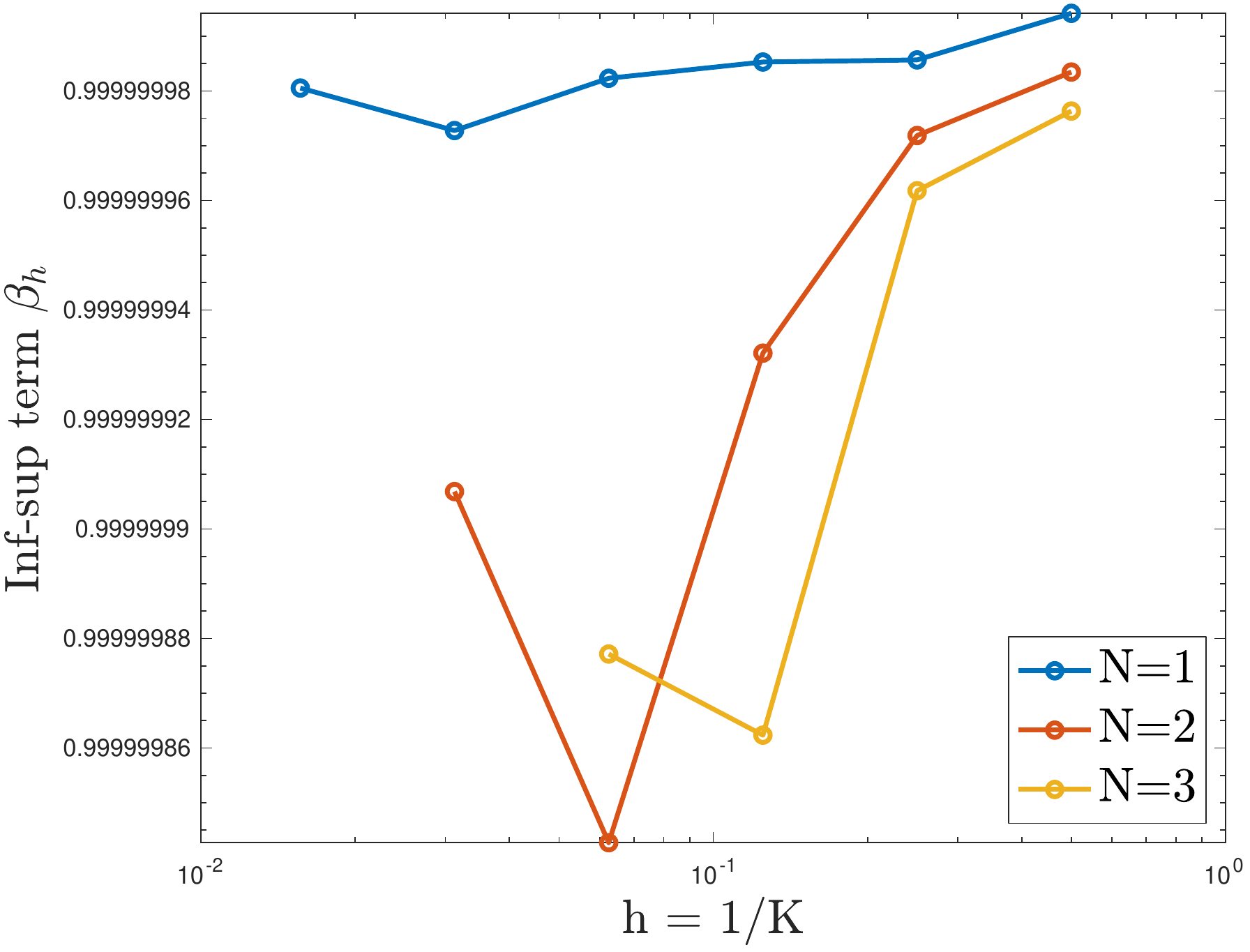}
\caption{Inf-sup term ${\beta}_h$ for $\Omega = [0,2]^2$ and polynomial degree $N=1,2,3$, with varying mesh refinement. Both the plots have same values. The plot on the right side is a zoom in of data.}
\label{fig:inf_sup_beta2}
\end{figure}
\section{Numerical tests}
In \figref{fig:inf_sup_beta1}, and \figref{fig:inf_sup_beta2}, we plot the value of the inf-sup term obtained using \eqref{eq:smallest_singular_value} for domain $\Omega = [0,1]^2$, and $\Omega = [0,2]^2$, respectively.
On the y-axis we have $\beta _h$, and on the x-axis we have the length of the element, $h=1/K$ for $K = 1 \;, 2  \;, 4  \;, 8  \;, 16  \;, 32 \;, 64$.
The plots on the right side in the figure are a zoom in view of the plots on the left side.
The numeric values for plots in \figref{fig:inf_sup_beta1} and \figref{fig:inf_sup_beta2} are given in \tabref{tab:smallest_singular_value1} and \tabref{tab:smallest_singular_value2} respectively.
In both the figures and the plots we observe that the numeric values of inf-sup term are all very close to $1.0$ and exact upto at least six decimal places, \VJ{which is in agreement with theoretical derivation}.
\begin{table}[!hbt]
	\caption{Numerical data for inf-sup term for $\Omega = [0,1]^2$.} \label{tab:smallest_singular_value1}
	\begin{tabular}{c|ccccccc}
		\hline
		 h    & $ N=1$        & $N=2$        & $N=3$       \\
		\hline
$1 / 2$	   	& 0.999999994172141 &   0.999999989168835 &   0.999999978662899 \\
$1 / 4$ 	&    0.999999987740597  &  0.999999978618959  &  0.999999947316136 \\
$1 / 8$	&    0.999999986310051   & 0.999999962318488  &  0.999999878165505 \\
$1 / 16$	&    0.999999980277492  &  0.999999916432749  &  0.999999891563428 \\
$1 / 32$	&    0.999999972413976  &  0.999999896705947      &               \\
$1 / 64$ &    0.999999980522211     &                              &      \\
		\hline
	\end{tabular}
\end{table}

\begin{table}[!hbt]
	\caption{Numerical data for inf-sup term for $\Omega = [0,2]^2$.} \label{tab:smallest_singular_value2}
	\begin{tabular}{c|cccccccc}
		\hline
		 h    & $ N=1$        & $N=2$        & $N=3$       \\
		\hline
$1 / 2$   & 0.999999994172141    & 0.999999983449168   & 0.999999976327674 \\
$1 / 4$   & 0.999999985661708   & 0.999999971854628   & 0.999999961774105 \\
$1 / 8$   & 0.999999985276638   & 0.999999932115646   & 0.999999862356188 \\
$1 / 16$   & 0.999999982313628   & 0.999999842706989  & 0.999999877166154 \\
$1 / 32$   & 0.999999972765685   & 0.999999906834163                   \\
$1 / 64$   & 0.999999980528361                                      \\
		\hline
	\end{tabular}
\end{table}

\section{Conclusions}
In this paper, we derive a theoretical estimate for the discrete inf-sup formulation and validate the value of the constant using finite dimensional spaces defined in \cite{2020jain}.
\VJ{The theoretical proof of the inf-sup term becomes straight forward when we use the appropriate norm on $K^\perp$ space, see \eqref{eq:norm_u1}.}
We evaluate the constant for \VJ{two} different test cases, i) a unit square domain $\Omega = [0,1]^2$, ii) a square domain $\Omega = [0,2]^2$.
It is shown that for all the cases the numerical value \VJ{is in agreement with the theoretical value}.

\bibliographystyle{elsarticle-num}
\bibliography{references}

\end{document}